\documentclass[11pt]{article}
\usepackage{stmaryrd}
\usepackage{tipa}
\usepackage{tikz}
\usepackage{amsmath}
\usepackage{amssymb}
\usepackage{pst-node}
\usepackage{color}
\usepackage{extarrows}

\usetikzlibrary{arrows,shapes,positioning}
\usetikzlibrary{decorations.markings}
\tikzstyle arrowstyle=[scale=1]
\tikzstyle directed=[postaction={decorate,decoration={markings, mark=at position 0.75 with {\arrow[arrowstyle]{stealth}}}}]
\tikzstyle redirected=[postaction={decorate,decoration={markings, mark=at position 0.35 with {\arrow[arrowstyle]{stealth}}}}]

\usepackage{authblk}
\usepackage{enumerate}
\usepackage{mathrsfs}
\usepackage{wrapfig}
\usepackage{float}

\begin{document}

\newtheorem{theorem}{Theorem}[section]
\newtheorem{corollary}[theorem]{Corollary}
\newtheorem{definition}[theorem]{Definition}
\newtheorem{conjecture}[theorem]{Conjecture}
\newtheorem{observation}[theorem]{Observation}
\newtheorem{problem}[theorem]{Problem}
\newtheorem{lemma}[theorem]{Lemma}
\newtheorem{newlemma}{New-Lemma}
\newtheorem{proposition}[theorem]{Proposition}
\newtheorem{construction}[theorem]{Construction}


\title{Some Early Results by Tutte Regarding the Cycle Double Cover Conjecture in 1948}
\author{Cun-Quan Zhang \\ Department of Mathematics, West Virginia University}
\date{}
\maketitle

\begin{abstract}OpenAI recently announced a proof of the Cycle Double Cover (CDC) Conjecture.
Most media reports have characterized it as a 50-year-old open problem.
In reality, according to a 1987 letter from Tutte to Fleischner, the Cycle Double Cover Problem
 has been open for at least 80 years.
Two early results regarding the CDC conjecture were established in one of Tutte's 1949 publications.

\noindent
{\bf Keywords:} Cycle Double Cover, Graph Theory, History of Mathematics.
\end{abstract}

\section{The CDC Conjecture Has Been Open for at Least 80 Years}
 OpenAI \cite{OpenAI2026} recently announced a proof of the Cycle Double Cover (CDC) Conjecture, which has long stood as one of the central and most challenging open problems in graph theory.
 However, most media reports have characterized it as a 50-year-old open problem.
  In reality, according to a 1987 letter from Bill Tutte to Herbert Fleischner, the Cycle Double Cover Problem
   is at least 80 years old.
  The original text of the letter reads as follows:

 \begin{quote}
\begin{flushright}
 ``July 22, 1987
\end{flushright}

 Dear Professor Fleischner:

\hspace{.3in}
 Thank you for your letter about the cycle double cover conjecture, I too have been puzzled to find an original reference.
 I think the conjecture is one that was well established in mathematical conversation long before anyone thought of publishing it.

\hspace{.3in}
 I don't remember referring to the conjecture in my own writings.
 The nearest I came was to show that bicubic graphs, and other cubic graphs with Tait colourings had circuit-double-covers, ``On the imbedding of linear graphs is surfaces'', Proc. London Math. Soc. Ser. 2, Vol. 51 (1948).

\begin{center}
All good wishes,

 (Signature)

 W. T. Tutte''
 \end{center}
 \end{quote}

\begin{figure}[H]
\begin{center}
\includegraphics[scale=0.82]{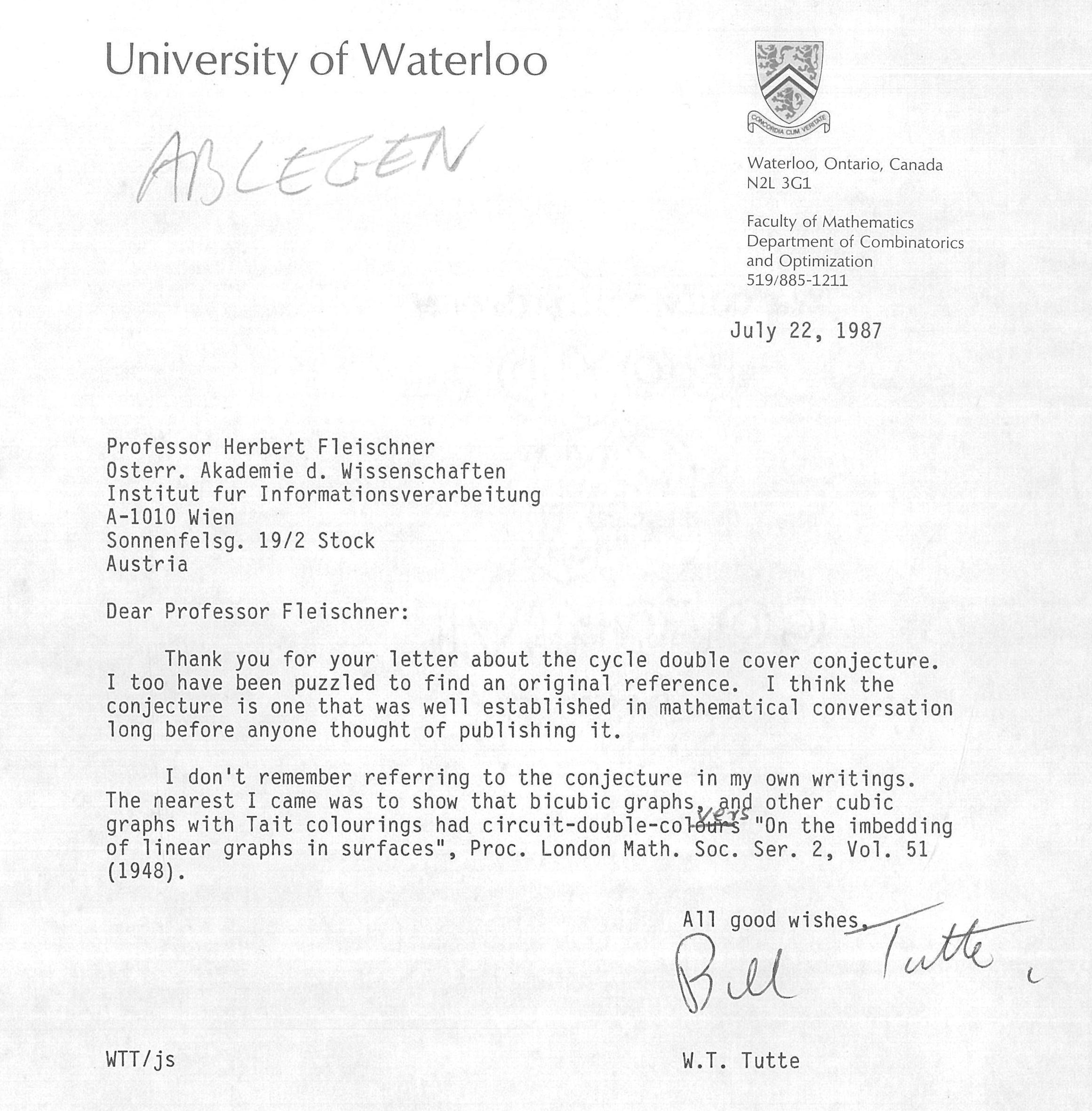}
\caption{\small\it A letter from Tutte to Fleischner on July 22, 1987}
\label{FIG: Tutte}
\end{center}
\end{figure}

 A photograph of the original letter is shown below (see Figure\ref{FIG: Tutte}).

 \section{Two Early Results Regarding CDC by Tutte in 1948}

 \begin{definition}
 An orientable \(k\)-even subgraph double cover of a graph \(G\) is a family \(\mathcal{F}\) of even subgraphs of \(G\) such that:

(1) \(\vert{}\mathcal{F}\vert{} = k\);

 (2) Every edge of \(G\) is covered by precisely two members of \(\mathcal{F}\);

 (3) Each member of \(\mathcal{F}\) has an Eulerian orientation and, furthermore, for each edge \(e\) of \(G\), \(e\) is contained in two members of \(\mathcal{F}\) with opposite directions.
     \end{definition}

     Below are two theorems established by Tutte in his {\em Proc. London Math. Soc.} paper:
     \begin{theorem}
     [\cite{Tutte1949}]
     Let \(G\) be a cubic graph. The graph \(G\) has an orientable \(3\)-even subgraph double cover if and only if \(G\) is bipartite.
     \end{theorem}

     \begin{theorem}[\cite{Tutte1949}]
     Let \(G\) be a cubic graph. The graph \(G\) has an orientable \(4\)-even subgraph double cover if and only if \(G\) is \(3\)-edge-colorable.\end{theorem}

By applying the concepts and notion
 of integer flow (introduced by Tutte  \cite{Tutte1949} \cite{Tutte1954}) and
  graph embedding, both theorems can be presented as follows.

\begin{corollary}[\cite{Tutte1949}]
  Let \(G\) be a cubic graph.
   Then the following statements are equivalent.

  (1) The graph \(G\) has an orientable \(3\)-even subgraph double cover;

  (2) \(G\) is bipartite;

  (3) $G$ admits a nowhere-zero $3$-flow;

  (4) $G$ has a strong embedding on an orientable surface ${\cal S}$ on which $G$ is $3$-face-colorable.
  \end{corollary}

  \begin{corollary}[\cite{Tutte1949}]
  Let \(G\) be a cubic graph.
   Then the following statements are equivalent.

  (1) The graph \(G\) has an orientable \(4\)-even subgraph double cover;

  (2) \(G\) is $3$-edge-colorable;

  (3) $G$ admits a nowhere-zero $4$-flow;

  (4) $G$ has a strong embedding on an orientable surface ${\cal S}$ on which $G$ is $4$-face-colorable.
  \end{corollary}

     \section{A Brief History Behind the Letter}
     The story traces back to a time before 2010, inside Fleischner's office in Vienna. Spread across his desk were several offprints and books, including the 1993 paper published by Alspach and myself in {\em Discrete Mathematics} \cite{Alspach1993}, its subsequent follow-up \cite{Alspach1994}, and my first book \cite{ZhangBook1}. It seemed as though they had been laid out specifically in anticipation of our conversation.
     Pointing to the pile of materials, he asked,
     ``Do you know how many years the CDC Conjecture has been open?''
      I paused for a moment before replying, ``It seems like it's been about 30 years?''
      He then opened a folder, carefully drew out a piece of stationery, and handed it to me.
      ``Take a look at this.''
      It was the original 1987 letter from Tutte mentioned above. \(\dots \)

      ``This material is incredibly precious. Thank you so much! I will definitely cite it in my next book,''
      I said. At that time, I was actively collecting materials for my second book \cite{ZhangBook2}. Fleischner immediately handed me a photocopy of the letter, saying, ``This is for you. Keep it.'' \(\dots \)

      \section{In Memory of Professor Herbert Fleischner}
      Prof. Fleischner passed away on October 7, 2025. He was long recognized as one of the leading figures in this field.
      Our connection began back in 1986, shortly after Alspach and I publicly announced the proof of the CDC Conjecture for Petersen-minor-free cubic graphs \cite{Alspach1993}.
      From that point forward, Fleischner and I remained in constant contact.
      Prior to the Covid-19 pandemic, I visited Vienna for what would be my fourth and final time.
      At that time, he was still actively preparing Volume 3 (or Part 2) of his seminal work, {\em ``Eulerian Graphs and Related Topics''} \cite{Fleischner1990, Fleischner1991}.
      He once mentioned to me that
      ``CDC will be one of the major subjects of the book.''
      We had all been eagerly anticipating this volume, which was sure to contain many of his original techniques as well as numerous open problems awaiting our resolution.

     \end{document}